\documentclass[a4paper,12pt]{amsart}
\usepackage{amssymb}
\usepackage{ifthen}
\nonstopmode
\numberwithin{equation}{section}
\newtheorem{thm}{Theorem}[section]
\newtheorem{cor}[thm]{Corollary}
\newtheorem{lem}[thm]{Lemma}
\newtheorem{prop}[thm]{Proposition}

\theoremstyle{remark}
\newtheorem{rem}[thm]{Remark}

\newcounter{alphabet}
\newcounter{tmp}

\newenvironment{pf}[1][]{%
 \vskip 3mm
 \noindent
 \ifthenelse{\equal{#1}{}}%
  {{\slshape Proof. }}%
  {{\slshape #1.} }%
 }%
{\qed\bigskip}

\begin{document}
\newcommand{\A}{{\mathcal A}}
\newcommand{\B}{{\mathcal B}}
\newcommand{\T}{{\mathcal T}}
\newcommand{\M}{{\mathcal M}}
\newcommand{\F}{{\mathcal F}}
\newcommand{\Om}{{\Omega}}
\newcommand{\es}{{\mathcal S}}
\newcommand{\R}{{\mathbb R}}
\newcommand{\C}{{\mathbb C}}
\newcommand{\K}{{\mathcal K}}
\newcommand{\E}{{\mathcal E}}
\newcommand{\eK}{{\bold K}}
\newcommand{\uhp}{{\mathbb H}}
\newcommand{\Z}{{\mathbb Z}}
\newcommand{\N}{{\mathbb N}}
\newcommand{\D}{{\mathbb D}}
\newcommand{\UCV}{{\mathcal{UCV}}}
\newcommand{\kUCV}{{k\text{-}\mathcal{UCV}}}
\newcommand{\Perron}{{\mathcal P}}
\newcommand{\sphere}{{\widehat{\mathbb C}}}
\newcommand{\image}{{\operatorname{Im}\,}}
\renewcommand{\Im}{{\operatorname{Im}\,}}
\newcommand{\Aut}{{\operatorname{Aut}\,}}
\newcommand{\real}{{\operatorname{Re}\,}}
\renewcommand{\Re}{{\operatorname{Re}\,}}
\newcommand{\kernel}{{\operatorname{Ker}\,}}
\newcommand{\id}{{\operatorname{id}}}
\newcommand{\mob}{{\text{\rm M\"{o}b}}}
\newcommand{\Int}{{\operatorname{Int}\,}}
\newcommand{\Ext}{{\operatorname{Ext}\,}}
\renewcommand{\mod}{{\operatorname{mod}}}
\newcommand{\stab}{{\operatorname{Stab}}}
\newcommand{\SL}{{\operatorname{SL}}}
\newcommand{\PSL}{{\operatorname{PSL}}}
\newcommand{\PSU}{{\operatorname{PSU}}}
\newcommand{\tr}{{\operatorname{tr}}}
\newcommand{\diam}{{\operatorname{diam}\,}}
\newcommand{\inv}{^{-1}}
\newcommand{\area}{{\operatorname{Area}\,}}
\newcommand{\eit}{{e^{i\theta}}}
\newcommand{\eint}{{e^{in\theta}}}
\newcommand{\emint}{{e^{-in\theta}}}
\newcommand{\dist}{{\operatorname{dist}}}
\newcommand{\arctanh}{{\operatorname{arctanh}}}
\newcommand{\const}{{\operatorname{const.}}}
\newcommand{\capa}{{\operatorname{Cap}}}
\newcommand{\hdim}{{\operatorname{H-dim}}}
\newcommand{\rad}{{\operatorname{rad}}}
\newcommand{\partialb}{{\partial_{\operatorname{b}}}}
\newcommand{\CD}{{\operatorname{CD}}}
\newcommand{\hm}{{\mathcal H}}
\newcommand{\hc}{{\mathcal L}}
\newcommand{\cube}{{\mathcal Q}}
\newcommand{\Log}{{\,\operatorname{Log}}}

\newcommand{\Der}{{\frak D}}
\newcommand{\X}{{\frak X}}
\newcommand{\Isom}{{\operatorname{Isom}}}
\newcommand{\Hom}{{\operatorname{Hom}}}
\newcommand{\der}{{\mathcal D}}
\renewcommand{\r}{{\varphi}}
\newcommand{\s}{{\psi}}
\newcommand{\z}{{\partial/\partial z}}
\newcommand{\zf}{{\frac\partial{\partial z}}}
\newcommand{\zb}{{\partial/\partial \bar z}}
\newcommand{\zbf}{{\frac\partial{\partial \bar z}}}
\newcommand{\mone}{{\mbox{-}1}}
\newcommand{\iu}{{\operatorname{i}}}
\newcommand{\gs}{{\Sigma}}
\newcommand{\poly}{{\mathcal P}}


\bibliographystyle{amsplain}
\title
[Invariant Schwarzian derivatives]
{
Invariant Schwarzian derivatives of higher order
}

\author[S.-A Kim]{Seong-A Kim}
\address{Department of Mathematics Education, Dongguk University \\
780-714, Korea}
\email{sakim@dongguk.ac.kr}  

\author[T. Sugawa]{Toshiyuki Sugawa}
\address{Graduate School of Information Sciences,
Tohoku University, Aoba-ku, Sendai 980-8579, Japan}
\email{sugawa@math.is.tohoku.ac.jp}

\begin{abstract}
We argue relations between the Aharonov invariants and Tamanoi's Schwarzian
derivatives of higher order and give a recursion formula for
Tamanoi's Schwarzians.
Then we propose a definition of invariant Schwarzian derivatives of
a nonconstant holomorphic map between Riemann surfaces with conformal metrics.
We show a recursion fomula also for our invariant Schwarzians.
\end{abstract}
\thanks{
The second author was supported in part by JSPS Grant-in-Aid for Scientific
Research (B), 17340039 and for Exploratory Research, 19654027.
}

\maketitle

\section{Introduction}

The Schwarzian derivative $S_f$
of a non-constant meromorphic function $f$ on a plane domain
is defined by
$$
S_f
=\frac{f'''}{f'}-\frac32\left(\frac{f''}{f'}\right)^2
=\left(\frac{f''}{f'}\right)'-\frac12\left(\frac{f''}{f'}\right)^2
=T_f'-\frac12T_f^2,
$$
where $T_f=f''/f'$ is the pre-Schwarzian derivative of $f.$
Note that $S_f(z)$ is holomorphic at $z=z_0$ when $f(z)$ is
locally univalent at $z=z_0,$ whereas $S_f(z)$ has a pole of order 2 at $z=z_0$
when $f(z)$ has a branch at $z=z_0.$

It is well recognized that the pre-Schwarzian derivative is crucially used
to construct a conformal mapping of the upper half-plane onto
a polygonal domain.
(This is the so-called Schwarz-Christoffel mapping.)
The Schwarzian derivative was introduced by Schwarz to
construct further a conformal mapping of the upper half-plane onto
a simply connected domain bounded by finitely many circular arcs.
After Nehari discovered univalence criteria of meromorphic
functions in terms of the Schwarzian derivative in the late 1940's,
Bers and Ahlfors found an intimate connection with quasiconformal mappings
and utilized it to embed Teichm\"uller spaces onto bounded domains
in complex Banach spaces (see, for example, \cite{Lehto:univ}).
Thus one may be tempted to define
higher-order analogues of the Schwarzian derivative.
Indeed, Aharonov \cite{Ahar69} and Tamanoi \cite{Tamanoi96} gave
definitions of higher-order analogues of the Schwarzian derivative.
Aharonov gave a necessary and sufficient condition for a nonconstant
meromorphic function on the unit disk to be univalent in terms
of his Schwarzians, whereas
Tamanoi studied combinatorial structures of his Schwarzians.
In Section 2, we briefly recall their definitions
and argue the relation between them.
We should also note here that Schippers \cite{Schip00} proposed yet another
definition of Schwarzian derivatives of higher order.
They fit the L\"owner theory and have nice properties.
However, as he noted in his paper, his Schwarzians have a different nature
from those of Aharonov and Tamanoi.
Thus we do not treat with Schipper's Schwarzians in this note.

One reason why the Schwarzian derivative is so useful is that
it has a nice invariance property.
On the other hand, Peschl and Minda introduced a sort of invariant
derivatives $D^nf$ of order $n$ for a holomorphic map $f$
between domains with conformal metrics
(see \cite{Schip07} and \cite{KS07diff}).
It is observed (\cite{KM93}, \cite{KM01}, \cite{MM97two}) that the quantity
$$
\gs f=\frac{D^3f}{D^1f} - \frac{3}{2}\left(\frac{D^2f}{D^1f}\right)^2,
$$
which has a form similar to $S_f,$ is essentially same as $S_f$
when the domains are standard ones (see also \cite{KS:univ}).
We generalize this idea to define invariant Schwarzian derivatives of
higher order in terms of the Peschl-Minda derivatives analogously.
One of our main results is a recursion formula of the invariant
Schwarzian derivatives of higher order (see Theorem \ref{thm:rec}).
These invariant Schwarzian derivatives and the recursion forumula
have applications to univalence criteria (see \cite{KS:univ} for details).

\section{Schwarzian derivatives of higher order}

Let $f$ be a nonconstant meromorphic function on a domain $D$ 
in the complex plane.
For $z\in D$ with $f(z)\ne\infty, f'(z)\ne0,$ we consider the quantity
$$
G(\zeta,z)=\frac{f'(z)}{f(\zeta)-f(z)}.
$$
We now expand it in the power series
$$
G(z+w,z)=\frac1w-\sum_{n=1}^\infty \psi_n[f](z)w^{n-1}
$$
for small enough $w.$
The quantities $\psi_n[f](z)$ were introduced by Aharonov \cite{Ahar69}
and called the {\it Aharonov invariants} by Harmelin \cite{Har82}.
Since the quantity
$$
\frac{\partial G}{\partial \zeta}(\zeta,z)
=-\frac{f'(z)f'(\zeta)}{(f(\zeta)-f(z))^2}
=\frac{-1}{(\zeta-z)^2}-\sum_{n=1}^\infty (n-1)\psi_n[f](z)(\zeta-z)^{n-2}
$$
is invariant under the M\"obius transformations of $f,$
we obtain $\psi_n[M\circ f]=\psi_n[f]$ for $n\ge2$ and 
a M\"obius transformation $M(z)=(az+b)/(cz+d).$
Thus these quantities can be defined even when $f(z)=\infty$
as long as $f$ is locally univalent at $z.$
Note that
$$
\psi_1[f](z)=\frac{f''(z)}{2f'(z)}
\quad\text{and}\quad
\psi_2[f](z)=\frac16\left[\frac{f'''(z)}{f'(z)}
-\frac32\left(\frac{f''(z)}{f'(z)}\right)^2\right].
$$
Thus $2!\psi_1[f]$ and $3!\psi_2[f]$ are the pre-Schwarzian derivative $T_f$
and the Schwarzian derivative $S_f$ of $f,$ respectively.
Thus, $\psi_n[f],~n=2,3,\dots,$ can be regarded as Schwarzian derivatives
of higher order.

Tamanoi \cite{Tamanoi96} proposed another definition of Schwarzian
derivatives of higher order.
Let $f$ be meromorphic in $D.$
Fix a point $z\in D$ with $f(z)\ne\infty, f'(z)\ne0,$
and take a M\"obius transformation $M_z$
so that
$$
M_z(0)=f(z),~
M_z'(0)=f'(z),~
M_z''(0)=f''(z).
$$
We will use later a concrete form of the inverse map $N_z=M_z\inv$ of $M_z:$
$$
N_z(t)=\frac{f'(z)(t-f(z))}{\frac12f''(z)(t-f(z))+f'(z)^2}.
$$
Then we expand the function $V=(M_z\inv\circ f)(z+w)=N_z(f(z+w))$ 
as a power series
\begin{equation}\label{eq:W}
V=\frac{f'(z)(f(z+w)-f(z))}{\frac12f''(z)(f(z+w)-f(z))+f'(z)^2}
=\sum_{n=0}^\infty S_n[f](z)\frac{w^{n+1}}{(n+1)!}
\end{equation}
around $w=0.$
The quantity $S_n[f]$ is called the Schwarzian derivative of
virtual order $n$ for $f$ (see \cite{Tamanoi96}).
By the choice of $M_z,$ we see that $S_0[f]=1, S_1[f]=0$ and $S_2[f]$ is the
classical Schwarzian derivative $S_f=f'''/f'-3(f''/f')^2/2.$
Also, by construction, $V$ and thus $S_n[f]$ are M\"obius invariant.
In particular, $S_n[f](z)$ can be defined even when $f(z)=\infty$
as long as $f$ is locally univalent at $z.$

Aharonov \cite{Ahar69} (see also \cite{Har82}) proved the recursion formula
$$
(n+1)\psi_n[f]=\psi_{n-1}[f]'+\sum_{k=2}^{n-2}\psi_k[f]\psi_{n-k}[f].
$$
We now show a similar formula for Tamanoi's Schwarzians.

\begin{prop}\label{prop:Sn}
$$
S_{n}[f]
=S_{n-1}[f]'+\tfrac12 S_2[f]\sum_{k=1}^{n-1} \binom{n}{k}S_{k-1}[f] S_{n-k-1}[f],
\quad n\ge 3.
$$
\end{prop}

\begin{pf}
We denote by $\dot N_z(t)$ the partial derivative of $N_z(t)$
with respect to $z,$ namely, $\dot N_z(t)=\partial_z N_z(t).$
The following formula is easily verified by a direct computation:
\begin{equation}\label{eq:N_w}
\dot N_z(t)=-1-\tfrac12 S_2[f](z) N_z(t)^2.
\end{equation}
We now compute partial derivatives of $V=N_z(f(z+w)):$
$$
\partial_w V=N_z'(f(z+w))f'(z+w),\quad
\partial_z V=N_z'(f(z+w))f'(z+w)+\dot N_z(f(z+w)).
$$
By \eqref{eq:N_w}, we have
\begin{equation}\label{eq:V}
\partial_z V-\partial_w V=\dot N_z(f(z+w))=-1-\tfrac12 S_2[f](z) V^2.
\end{equation}
Compare with the similar formula (2.8) in \cite{Ahar69}.
We now substitute \eqref{eq:W} into the last formula to obtain
\begin{align*}
&~\sum_{n=0}^\infty S_n[f]'(z)\frac{w^{n+1}}{(n+1)!}
-\sum_{n=0}^\infty S_n[f](z)\frac{w^{n}}{n!} \\
=&~-1+\sum_{n=1}^\infty \big\{ S_{n-1}[f]'(z)-S_{n}f(z)\big\}\frac{w^{n}}{n!} \\
=&~-1-\tfrac12 S_2[f](z)]\sum_{k=1}^\infty\sum_{l=1}^\infty 
S_{k-1}[f](z)S_{l-1}[f](z)\frac{w^{k+l}}{k!\,l!} \\
=&~-1-\tfrac12 S_2[f](z)\sum_{n=2}^\infty\sum_{k=1}^{n-1} 
\binom{n}{k}S_{k-1}[f](z)S_{n-k-1}[f](z)\frac{w^{n}}{n!}.
\end{align*}
By comparing the coefficients of $w^n,$ we obtain the required relation.
\end{pf}

In particular, we have the following result.
Here and hereafter, $\Z$ denotes the ring of integers.
Note that a similar result was obtained by Tamanoi
(see \cite[Theorem 6-4]{Tamanoi96}).

\begin{cor}
$S_{n}[f]$ is expressed in the form 
$S_{n-1}[f]'+P(S_2[f],\dots, S_{n-4}[f],S_{n-2}[f])$ 
for $n\ge3,$ where $P(x_2,x_3, \dots, x_{n-5}, x_{n-4}, x_{n-2})$ is a polynomial
in $x_2, x_3, \dots, x_{n-5}, x_{n-4}, x_{n-2}$ 
with non-negative coefficients in $\Z.$
\end{cor}

\begin{pf}
When $n$ is odd, we can write
$$
S_{n}[f]=S_{n-1}[f]'+S_2[f]\sum_{k=1}^{\frac{n-1}2} 
\binom{n}{k}S_{k-1}[f] S_{n-k-1}[f].
$$
When $n$ is even, we have
$$
S_{n}[f]=S_{n-1}[f]'+S_2[f]\sum_{k=1}^{\frac n2-1} 
\binom{n}{k}S_{k-1}[f] S_{n-k-1}[f]
+\frac12\binom{n}{\frac n2}S_2[f]S_{\frac n2-1}[f]^2.
$$
Since
$$
2^n=(1+1)^n=\sum_{k=0}^n\binom{n}{k}=
2\sum_{k=0}^{\frac n2-1} \binom{n}{k}+\binom{n}{\frac n2},
$$
we see that $\binom{n}{n/2}$ is an even number for $n\ge3.$
Since $S_1[f]=0,$ the assertion follows.
\end{pf}

Let us write down first several nontrivial Schwarzians:
\begin{align*}
S_3[f]&= S_2[f]' \\
S_4[f]&= S_3[f]'+4S_2[f]^2 \\
S_5[f]&= S_4[f]'+5S_2[f]S_3[f] \\
S_6[f]&= S_5[f]'+6S_2[f]S_4[f]+10S_2[f]^3.
\end{align*}

It should be noted here that $S_n[f]=0$ does not necessarily imply 
$S_{n+1}[f]=0.$
For example, consider the function $f(z)=e^{az}$ for a constant $a\ne0.$
Then $S_2[f]=-a^2/2.$
Therefore, $S_3[f]=0$ but $S_4[f]=a^4.$

We conclude the present section with the relationship between
the Aharonov invariants and Tamanoi's Schwarzian derivatives.
For convenience, we set 
$$
\sigma_n[f]=\frac{S_n[f]}{(n+1)!}\quad\text{for}~n=0,1,2,\dots.
$$
By the elementary formula
$$
(G+\psi_1[f])V=\left(1-\sum_{n=2}^\infty\psi_n[f](z)w^n\right)
\sum_{n=0}^\infty\sigma_n[f](z)w^n
=1,
$$
we obtain the following.

\begin{prop}
$$
\sigma_n[f]=\psi_n[f]+\sum_{k=2}^{n-2}\psi_k[f]\sigma_{n-k}[f],
\quad n\ge 2.
$$
\end{prop}

\begin{cor}
$\sigma_n[f]$ can be expressed as a polynomial of
$\psi_2[f],\dots,\psi_n[f]$ with non-negative coefficients in $\Z.$
\end{cor}

For example, we have
\begin{align*}
\sigma_2[f]&=\psi_2[f], \\
\sigma_3[f]&=\psi_3[f], \\
\sigma_4[f]&=\psi_4[f]+\psi_2[f]^2, \\
\sigma_5[f]&=\psi_5[f]+2\psi_2[f]\psi_3[f] \\
\sigma_6[f]&=\psi_6[f]+\psi_2[f]^3+\psi_3[f]^2+2\psi_2[f]\psi_4[f].
\end{align*}

\section{Expression of $S_n[f]$ in terms of the quotients $f^{(k)}/f'$}

Let $R=\Z[\tfrac12]$ be the ring generated by $1/2$ over $\Z.$
We consider the ring $R[x_1,x_2,\dots]$ of polynomials
of infinitely many indeterminates $x_1,x_2,\dots$ over $R.$
The weight of a monomial $x_{j_1}\cdots x_{j_k}$ is defined to be
the number $j_1+\cdots+j_k.$
Let $\poly_m$ be the sub $R$-module of $R[x_1,x_2,\dots]$
generated by monomials of  weight $m.$
A polynomial $P\in R[x_1,x_2,\dots]$ is said to be {\it of weight} $m$
if $P\in\poly_m.$
It is easy to see that $\poly=\sum_{m=0}^\infty \poly_m$ becomes a graded ring.
We denote by $\E$ the $R$-derivation on $R[x_1,x_2,\dots]$ defined by
$$
\E P=\sum_{n=1}^\infty (x_{n+1}-x_1x_n)\frac{\partial P}{\partial x_n}.
$$
Note that $\E$ maps $\poly_m$ into $\poly_{m+1}.$
We now define polynomials $P_n\in R[x_1,x_2,\dots],~n=0,1,2,\dots,$ 
inductively by $P_0=1, P_1=0, P_2=x_2-(3/2)x_1^2$ and
\begin{equation}\label{eq:Pn}
P_{n}=\E P_{n-1}
+\frac12 P_2\sum_{k=1}^{n-1} \binom{n}{k}P_{k-1}P_{n-k-1},\quad n\ge3.
\end{equation}
For instance,
\begin{align*}
P_3&=x_3-4x_1x_2+3x_1^3, \\
P_4&=x_4-5x_1x_3+5x_1^2x_2, \\
P_5&=x_5-6x_1x_4+\frac{15}2 x_1^2x_3-10x_1x_2^2+30x_1^3x_2-\frac{45}2x_1^5, \\
P_6&=x_6-7x_1x_5+\frac{21}2 x_1^2x_4-35x_1x_2x_3+\frac{105}2x_1^3x_3
+105x_1^2x_2^2-210x_1^4x_2-\frac{315}4x_1^6, \\
P_7&=x_7-8x_1x_6+14x_1^2x_5-56x_1x_2x_4+84x_1^3x_4-35x_1x_3^2+420x_1^2x_2x_3 \\
&\quad -420x_1^4x_3-420x_1^3x_2^2+420x_1^5x_2.
\end{align*}

\begin{lem}\label{lem:Pn}
The above polynomials $P_n,~n=0,1,2,\dots,$ satisfy the following properties:
\begin{enumerate}
\item[(i)]
$P_n$ is of weight $n.$
\item[(ii)]
$P_n\in R[x_1,\dots, x_{n}].$
\item[(iii)]
$\displaystyle\sum_{k=1}^{n} kx_k\frac{\partial P_n}{\partial x_k}=nP_n.$
\end{enumerate}
\end{lem}

\begin{pf}
Property (i) can easily be checked by induction on $n.$
Property (ii) follows from (i).
We thus prove only property (iii).

For a monomial $A=x_1^{e_1}\dots x_{n}^{e_{n}},$ we compute
$$
x_k\frac{\partial A}{\partial x_k}=e_k A.
$$
Therefore,
$$
\sum_{k=1}^{n} jx_j\frac{\partial A}{\partial x_k}=
(e_1+2e_2+\dots +n e_{n})A.
$$
Note here that the weight of $A$ is
$e_1+2e_2+\dots +ne_{n}.$
Since $P_n$ is given as a linear combination of
monomials of weight $n,$
the assertion is now clear.
\end{pf}

By property (ii) above, we may think of $P_n$ as a function
$P_n(x_1,\dots,x_{n})$ of $x_1,\dots,x_{n}.$
Let
$$
q_n[f]=\frac{f^{(n+1)}}{f'}, \quad n=1,2,\dots.
$$
Then the principal result in this section is the following.

\begin{thm}\label{thm:Pn}
$$
S_n[f]=P_n(q_1[f],q_2[f],\dots,q_n[f]),\quad n\ge0.
$$
\end{thm}

\begin{pf}
For $n=0,1,2,$ this is clear by definition.
If this is true for $1,2,\dots, n-1,$ then
$$
S_{n-1}[f]'=P_{n-1}(q_1[f],\dots,q_{n-1}[f])'
=\sum_{k=1}^{n-1}\frac{\partial P_{n-1}}{\partial x_k}(q_1[f],\dots,q_{n-1}[f])
q_k[f]'.
$$
Since
\begin{equation}\label{eq:q}
q_k[f]'=q_{k+1}[f]-q_1[f]q_k[f],
\end{equation}
the relation $S_{n-1}[f]'=\E P_{n-1}(q_1[f],\dots,q_n[f])$ holds.
Now we use Proposition \ref{prop:Sn} to show the assertion by induction.
\end{pf}

\section{Invariant Schwarzian derivatives}

In this section, we first recall the definition of a sort of
invariant derivatives for holomorphic maps
between plane domains (or, more generally, Riemann surfaces)
with (smooth) conformal metrics.
These were introduced by Peschl \cite{Peschl55} when the domains
are either the unit disk, the complex plane or 
the Riemann sphere with canonical metrics.
Later the notion was generalized by Minda for general conformal metrics.
We call those derivatives the Peschl-Minda derivatives
and detailed accounts were recently supplied by Schippers \cite{Schip07}
and the authors \cite{KS07diff}.

For simplicity, we consider only plane domains in the present note.
However, the notions below can easily be extended for Riemann surfaces
as we will make a remark on it later.

Let $\Omega$ and $\Omega'$ be plane domains with (smooth) conformal metrics
$\rho=\rho(z)|dz|$ and $\sigma=\sigma(w)|dw|,$ respectively.
We first define the $\rho$-derivative of
a smooth function $\r$ on $\Omega$ by
$$
\partial_\rho\r
=\frac1{\rho(z)}\frac{\partial\r(z)}{\partial z}.
$$
For a holomorphic map $f:\Omega\to\Omega',$ we define invariant differential
operators
$D^nf=D^n_{\sigma,\rho}f$ of order $n$ with respect to $\rho$ and $\sigma$
inductively by
\begin{align*}
D^1_{\sigma,\rho}f&=\frac{\sigma\circ f}{\rho}f' \\
D^{n+1}_{\sigma,\rho}f&
=\left[\partial_\rho-n\partial_\rho(\log\rho)
+(\partial_\sigma\log\sigma)\circ f\cdot 
D^1_{\sigma,\rho}f\right]D^n_{\sigma,\rho}f \quad(n\ge1).
\end{align*}
See \cite{KS07diff} or \cite{Schip07} for details.

The quotient $Q_f=D^2f/D^1f$ with variable metrics
was effectively used by Ma, Minda and others
in the geometric study of analytic maps between plane domains (see, for instance,
\cite{KM01} and \cite{KS04HC} and references therein).
We also cosinder its higher-order analogues:
$$
Q^nf=\frac{D^{n+1}f}{D^1f},\quad n\ge1.
$$
When we need to indicate the metrics, we write $Q_{\sigma,\rho}^nf$
instead of $Q^nf.$

Let $P_n(x_1,\dots,x_n)$ be the polynomial defined in the previous section.
We define the {\it invariant Schwarzian derivative} $\gs^nf$
of virtual order $n$ for $f$ by
$$
\gs^nf=P_n(Q^1f,\dots,Q^nf).
$$
To indicate the metrics involved, we sometimes write 
$\gs^nf=\gs_{\rho,\sigma}^nf.$
Note that $\gs^nf$ reduces to $S^n[f]$ when $\rho=\sigma=|dz|.$
We have the following invariance property for these quantities.

\begin{lem}\label{lem:inv}
Let $\Omega, \hat\Omega, \Omega', \hat\Omega'$ 
be plane domains with smooth conformal
metrics $\rho, \hat\rho, \sigma, \hat\sigma,$ respectively.
Suppose that locally isometric holomorphic maps
$g:\hat\Omega\to \Omega$ and $h:\Omega'\to\hat\Omega'$ are given.
Then, for a non-constant holomorphic map $f:\Omega\to\Omega',$ the formulae
\begin{align*}
Q^n_{\hat\sigma,\hat\rho}(h\circ f\circ g)
&=(Q^n_{\sigma,\rho}f)\circ g
\cdot\left(\frac{g'}{|g'|}\right)^n \\
\gs^n_{\hat\sigma,\hat\rho}(h\circ f\circ g)
&=(\gs^n_{\sigma,\rho}f)\circ g
\cdot\left(\frac{g'}{|g'|}\right)^n
\end{align*}
are valid on $\hat\Omega.$
\end{lem}

\begin{pf}
By \cite[Lemma 3.6]{KS07diff}, we have
$$
D^n_{\hat\sigma,\hat\rho}(h\circ f\circ g)
=\left(\frac{h'}{|h'|}\right)\circ f\circ g\cdot(D^n_{\sigma,\rho}f)\circ g
\cdot\left(\frac{g'}{|g'|}\right)^n.
$$
Thus the assertion for $Q^n$ follows immediately.
To prove that for $\gs^n,$ it is enough to observe the identity
$$
P_n(\alpha x_1,\alpha^2 x_2,\dots, \alpha^n x_n)
=\alpha^n P_n(x_1,x_2,\dots,x_n)
$$
for $\alpha\in\C,$ which can be seen easily.
\end{pf}

\begin{rem}
By the above lemma, the quantities $Q^nf$ and $\gs^nf$
can be defined as $(\frac n2,-\frac n2)$-forms on the Riemann surface $R$
for a non-constant holomorphic map $f:R\to R'$ between
Riemann surfaces $R$ and $R'$ with conformal metrics.
In particular, $|Q^nf|$ and $|\gs^nf|$ are independent of
the particular choices of local coordinates and thus can be regarded
as functions on $R.$
\end{rem}

The next result is an analogue of \eqref{eq:q}.

\begin{lem}\label{lem:dQ}
$$
\partial_\rho(Q^nf)=Q^{n+1}f-\big[Q^1f-n\partial_\rho\log\rho\big] Q^nf.
$$
\end{lem}

\begin{pf}
Recall that
$$
D^{n+1}f=\partial_\rho(D^nf)+\big[-n\partial_\rho\log\rho+
(\partial_\sigma\log\sigma)\circ f\cdot D^1f\big]D^nf
$$
for $n\ge1.$
By dividing both sides by $D^1f,$ we have
$$
Q^{n}f=\frac{\partial_\rho(D^nf)}{D^1f}
+\big[-n\partial_\rho\log\rho
+(\partial_\sigma\log\sigma)\circ f\cdot D^1f\big] Q^{n-1}f.
$$
Since
\begin{align*}
\frac{\partial_\rho(D^nf)}{D^1f}
&=\frac{\partial_\rho(Q^{n-1}f\cdot D^1f)}{D^1f} \\
&=\partial_\rho(Q^{n-1}f)+Q^{n-1}f\cdot\frac{\partial_\rho D^1f}{D^1f} \\
&=\partial_\rho(Q^{n-1}f)+Q^{n-1}f\big[Q^1f+\partial_\rho\log\rho
-(\partial_\sigma\log\sigma)\circ f\cdot D^1f\big],
\end{align*}
we obtain the assertion for $n-1.$
\end{pf}

We are now able to show the following result, which is a generalization
of Proposition \ref{prop:Sn}.

\begin{thm}\label{thm:rec}
Let $f$ be a non-constant holomorphic map between plane domains
$\Omega$ and $\Omega'$ with conformal metrics $\rho$ and $\sigma,$
respectively. Then
$$
\gs^{n}f=\big(\partial_\rho-(n-1)\partial_\rho\log\rho\big)\gs^{n-1}f
+\frac12\gs^2f\sum_{k=1}^{n-1} \binom{n}{k}\gs^{k-1}f \gs^{n-k-1}f,
\quad n\ge 3.
$$
\end{thm}

\begin{pf}
By definition, we compute
$$
\partial_\rho\gs^{n-1}f
=\partial_\rho P_{n-1}(Q^1f,\dots,Q^{n-1}f)
=\sum_{k=1}^{n-1}\frac{\partial P_{n-1}}{\partial x_k}(Q^1f,\dots,Q^{n-1}f)
\cdot\partial_\rho Q^{k}f.
$$
We now substitute the relation in Lemma \ref{lem:dQ} into the above
to get
\begin{align*}
\partial_\rho\gs^{n-1}f
&=\sum_{k=1}^{n-1}\frac{\partial P_{n-1}}{\partial x_k}(Q^1f,\dots,Q^{n-1}f)
\big[Q^{k+1}f-\big[Q^1f-k\partial_\rho\log\rho\big] Q^{k}f\big] \\
&=\sum_{k=1}^{n-1}\frac{\partial P_{n-1}}{\partial x_k}(Q^1f,\dots,Q^{n-1}f)
\big[Q^{k+1}f-Q^1fQ^{k}f\big] \\
&\quad +\partial_\rho\log\rho
\sum_{k=1}^{n-1}k Q^{k}f\cdot
\frac{\partial P_{n-1}}{\partial x_k}(Q^1f,\dots,Q^{n-1}f) \\
&=(\E P_{n-1})(Q^1f,\dots,Q^{n}f)
+(n-1)\partial_\rho\log\rho\cdot P_{n-1}(Q^1f,\dots,Q^{n-1}f),
\end{align*}
where we used the definition of $\E$ and 
property (iii) in Lemma \ref{lem:Pn}.
We finally recall the defining relation \eqref{eq:Pn} of $P_n$ to obtain
\begin{align*}
\partial_\rho\gs^{n-1}f
&=\left\{P_{n}
-\frac12 P_2\sum_{k=1}^{n-1} \binom{n}{k}
P_{k-1}P_{n-k-1}\right\}(Q^1f,\dots,Q^{n}f) \\
&\quad +(n-1)\partial_\rho\log\rho\cdot P_{n-1}(Q^1f,\dots,Q^{n-1}f) \\
&=\gs^{n}f-\frac12 \gs^2f\sum_{k=1}^{n-1}\binom{n}{k} \gs^{k-1}f\gs^{n-k-1}f
+(n-1)\partial_\rho\log\rho\cdot \gs^{n-1}f.
\end{align*}
Thus the assertion has been shown.
\end{pf}

We remark that the first author \cite{Kim:gs} gives a generating function
of the invariant Schwarzian derivatives $\gs^nf$ and proves Theorem \ref{thm:rec}
based on a relation similar to \eqref{eq:V} for $f:\C_\delta\to\C_\varepsilon.$
Here $\delta, \varepsilon=+1, 0,$ or $-1,$ and $\C_{+1}=\sphere,
\C_0=\C$ and $\C_{-1}=\D$ and the domain $\C_\delta$ is equipped with
the standard metric $\lambda_\delta=|dz|/(1+\delta|z|^2)$ for $\delta=+1,0,-1.$

\def\cprime{$'$} \def\cprime{$'$} \def\cprime{$'$}
\providecommand{\bysame}{\leavevmode\hbox to3em{\hrulefill}\thinspace}
\providecommand{\MR}{\relax\ifhmode\unskip\space\fi MR }
\providecommand{\MRhref}[2]{%
  \href{http://www.ams.org/mathscinet-getitem?mr=#1}{#2}
}
\providecommand{\href}[2]{#2}

\end{document}